\documentclass[12pt]{amsart}
\usepackage{amssymb}
\input amssym.def
\usepackage{amsmath, amsfonts,hyperref}
\usepackage{amscd}
\usepackage[mathscr]{eucal}
\newcommand{\Z}{{\mathbb Z}}

\newcommand{\starsum}{\mathop{\enspace{\sum}^{\ast}}}

\newtheorem{thm}{Theorem}
\newtheorem{lem}{Lemma}

\newtheorem{rmk}{Remark}
\newtheorem{defn}{Definition}
\newcommand{\thmref}[1]{Theorem~\ref{#1}}

\newcommand{\lemref}[1]{Lemma~\ref{#1}}

\begin{document}

\title[Finite Ramanujan expansions and shifted convolution sums]
{Finite Ramanujan expansions and shifted convolution sums of
arithmetical functions}

\author{Giovanni Coppola, M. Ram Murty and Biswajyoti Saha}

\address{Giovanni Coppola\\ \newline
Universit\'a degli Studi di Napoli,
Complesso di Monte S. Angelo-Via Cinthia
80126 Napoli (NA), Italy}
\email{giovanni.coppola@unina.it}

\address{M. Ram Murty\\ \newline
Department of Mathematics, Queen's University,
Kingston, Ontario, K7L 3N6, Canada}
\email{murty@mast.queensu.ca}
\thanks{Research of the second author was partially supported by an
NSERC Discovery grant.}

\address{Biswajyoti Saha\\ \newline
School of Mathematics, Tata Institute of Fundamental Research,
Homi Bhabha Road, Navy Nagar, Mumbai, 400 005, India}
\email{biswa@math.tifr.res.in}



\date{\today}

\begin{abstract}
For two arithmetical functions $f$ and $g$, we study the convolution sum
of the form $\sum_{n \le N} f(n) g(n+h)$ in the context of its asymptotic formula
with explicit error terms. Here we introduce the concept of finite
Ramanujan expansion of an arithmetical function and extend our earlier
works in this setup.
\end{abstract}

\maketitle

\section{Introduction}

In 1918, Ramanujan studied \cite{SR} the following
sum of roots of unity.

\begin{defn}
For positive integers $r,n$,
$$
 c_r(n):=\sum_{a\in (\Z/r\Z)^*}\zeta_r^{an},
$$
where $\zeta_r$ denotes a primitive $r$-th root of unity.
\end{defn}

These sums are now known as Ramanujan sums.
It is also possible to write $c_r(n)$ in
terms of the M\"obius function $\mu$ (see \cite{RM}). One has
\begin{equation}\label{crn-exp1}
c_r(n)=\sum_{d|n, d|r} \mu(r/d) d
\end{equation}
for any positive integers $r,n$.

Ramanujan studied these sums in the context of
point-wise convergent series expansions of the form
$\sum_r a_rc_r(n)$ for various arithmetical functions.
Such expansions are now known as Ramanujan expansions. More precisely:

\begin{defn}
We say an arithmetical function $f$
admits a Ramanujan expansion  (in the sense of Ramanujan) if for each $n$,
$f(n)$ can be written as a convergent series
of the form
$$
f(n)= \sum_{r \ge 1} \hat f(r)c_r(n)
$$
for appropriate complex numbers $\hat f(r)$.
The number $\hat f(r)$ is said to be the $r$-th Ramanujan coefficient
of $f$ with respect to this expansion.
\end{defn}

Shifted convolution sums are ubiquitous in number theory and
recently such sums have been studied for functions with
absolutely convergent Ramanujan expansions. It
has been done systematically in \cite{GMP,MS,CMS,BS2}.
For two arithmetical functions $f$ and $g$ we study the convolution
sum of the form $\sum_{n \le N} f(n) g(n+h)$. In this article we
introduce the concept of finite Ramanujan
expansion of an arithmetical function. This idea particularly
enables us to avoid technical infinite sums and obtain an
asymptotic formula with explicit error terms for the convolution sum
$C_{f,g}(h):=\sum_{n \le N} f(n) g(n+h)$ for some fixed positive integer $N$
and non-negative integer $h$.

Let us write our functions $f$ and $g$ as
$$
f(n)=\sum_{d|n}f'(d) \ \ \text{and} \ \
g(n)=\sum_{d|n}g'(d),
$$
where $f' := f \ast \mu$ and $g' := g \ast \mu$.
Here $\mu$ denotes the  M\"obius function and
$\ast$ denotes the Dirichlet convolution. Then
$$
C_{f,g}(h)=\sum_{n\le N} \sum_{d|n}f'(d) \sum_{q|n+h}g'(q)
=\sum_{d\le N}f'(d)\sum_{q\le N+h}g'(q)
\sum_{{n\le N}\atop {{n\equiv 0\bmod d}\atop {n+h\equiv 0\bmod q}}}1.
$$
Hence from the point of view of studying the convolution sums, we may put
\begin{equation}\label{f}
f(n)=\sum_{d|n, d\le N}f'(d)
\end{equation}
and
\begin{equation}\label{g}
g(n+h)=\sum_{d|n+h, d\le N+h}g'(d),
\end{equation}
i.e. we are enforcing that $f'(n)$ vanishes if $n>N$ and $g'(n)$ vanishes
if $n>N+h$. This will definitely change the values our $f(n),g(n)$ that we
started with, but only for $n>N$ and $n>N+h$ respectively. Hence
this will not alter our convolution sum $\sum_{n \le N} f(n) g(n+h)$.
At this point we note the following interesting property satisfied by
Ramanujan sums.

\begin{lem}\label{l1}
$$
\frac{1}{d} \sum_{r|d} c_r(n)= \begin{cases}
1 & \text{if } \ d|n,\\
0 & \text{otherwise}.
\end{cases}
$$
\end{lem}

For a proof see Section \ref{lemmas}. Now, using \lemref{l1}, we get
$$
f(n) = \sum_{d|n, d\le N}f'(d) =  \sum_{d\le N} f'(d) \frac{1}{d} \sum_{r|d} c_r(n)
=  \sum_{r \le N} c_r(n) \left( \sum_{r|d, d\le N} \frac{f'(d) }{d} \right)
=  \sum_{r \le N} \hat f(r) c_r(n),
$$
where
\begin{equation}\label{hat f}
\hat f(r) := \sum_{r|d, d\le N} \frac{f'(d) }{d}.
\end{equation}
Similarly
$$
g(n+h) = \sum_{s \le N+h} \hat g(s) c_s(n+h),
$$
where
\begin{equation}\label{hat g}
\hat g(s) := \sum_{s|d, d\le N+h} \frac{g'(d) }{d}.
\end{equation}

Thus we obtain a finite series expansion for our functions $f$ and $g$
which is more like a Ramanujan expansion. This we refer to as finite Ramanujan
expansion relative to $N$ and $h$. Note that such kind of an expansion depends on the
fixed parameters $N$ and $h$. This helps us to avoid dealing with
infinite sums, which was not the case in \cite{MS,CMS,BS2}. From now on, all these
above notations will be used freely without referring to them.

Using the dual M\"obius inversion
formula (see page 4 of \cite{CM}) it is also possible to express
$f'$ in terms of $\hat f$. We have
\begin{equation}\label{f'}
f'(r) = r \sum_{r|d, d\le N} \mu(d/r) \hat f(d)
\end{equation}
and
\begin{equation}\label{g'}
g'(s) = s \sum_{s|d, d\le N+h} \mu(d/s) \hat g(d).
\end{equation}

For arithmetical functions with usual Ramanujan expansion that
are absolutely convergent, the following theorem was proved in \cite{CMS}.

\begin{thm}[Coppola-Murty-Saha]\label{CMS}
 Suppose that $f$ and $g$ are two arithmetical functions with absolutely
 convergent Ramanujan expansions (in the sense of Ramanujan):
 $$
 f (n) = \sum_{r \ge 1} \hat f(r) c_r(n), \phantom{mm}
 g(n) = \sum_{s \ge 1} \hat g(s)c_s(n)
 $$
 respectively. Further suppose that
 $$
 \big|\hat f(r)\big|,\big|\hat g(r)\big| \ll \frac{1}{r^{1 +\delta}}
 $$
 for some $\delta > 0$ and $h$ is a non-negative integer.
 Then we have
  $$
 \sum_{n \le N} f(n) g(n+h) =
 \begin{cases}
N \sum_{r \ge 1} \hat f(r) \hat g(r) c_r(h) +
O(N^{1-\delta}(\log N)^{4-2\delta}) & ~ \mbox{ if } ~ \delta < 1,\\
N \sum_{r \ge 1} \hat f(r) \hat g(r) c_r(h) +
O(\log^3 N) & ~ \mbox{ if } ~ \delta = 1,\\
N \sum_{r \ge 1} \hat f(r) \hat g(r) c_r(h) +
O(1) & ~ \mbox{ if } ~ \delta > 1.
 \end{cases}
 $$
\end{thm}

Now suppose we impose the conditions 
\begin{equation}\label{hypo-1}
\big|\hat f(r)\big|,\big|\hat g(r)\big| \ll \frac{1}{r^{1 +\delta}}
\ \ \text{for some} \ \ \delta>0
\end{equation}
on the coefficients of the finite Ramanujan expansions of the arithmetic
functions $f$ and $g$, defined in \eqref{hat f} and \eqref{hat g}. Using
the dual M\"obius inversion formula, these conditions can be rewritten
equivalently as
\begin{equation}\label{hypo-2}
\big|f'(r)\big|,\big|g'(r)\big| \ll \frac{1}{r^{\delta}}
\ \ \text{for some} \ \ \delta>0.
\end{equation}
With these conditions in place we can derive a theorem that is analogous
to \thmref{CMS} and we are also able to improve the error term in the case
of $\delta \le 1$ by certain exponents of $\log N$.

\begin{thm}\label{finite}
Let $N$ be a positive integer and $f$ and $g$ be two arithmetical functions
for which we want to estimate the shifted convolution sums
$$
\sum_{n \le N} f(n) g(n+h)
$$
for a positive integer $h$. Further suppose that
$$
\big|\hat f(r)\big|,\big|\hat g(r)\big| \ll \frac{1}{r^{1 +\delta}}
\ \ \text{for some} \ \  \delta > 0,
$$
where $\hat f(r),\hat g(r)$ are as in \eqref{hat f}
and \eqref{hat g}. Then we have
$$
\sum_{n \le N} f(n) g(n+h) = N \sum_{r=1}^{\infty}\hat{f}(r)\hat{g}(r)c_r(h)
+ O_{\delta,h}\Big(N^{1-\delta} \log^2 N+1\Big).
$$
\end{thm}

The study of shifted convolution sums in the context of arithmetical functions
with absolutely convergent Ramanujan expansions was initiated by Gadiyar, Murty
and Padma in \cite{GMP}. The authors in \cite{GP} showed that if we ignore
convergence questions, a Ramanujan expansion of the function
$\frac{\phi(n)}{n}\Lambda(n)$, which is due to Hardy, can be used to derive
the Hardy-Littlewood conjecture about prime tuples. Later in \cite{GMP}
it was investigated whether the work in \cite{GP} can be justified
for arithmetical functions with absolutely convergent Ramanujan expansion,
under certain hypothesis on the Ramanujan coefficients. The objective of \cite{BS2}
was to reach the minimality of such hypothesis. In that quest the last author
proved the following theorem in \cite{BS2}.

\begin{thm}[Saha]\label{BS}
Let $f,g$ be two arithmetical functions with absolutely
convergent Ramanujan expansions (in the sense of Ramanujan)
$$
f (n) = \sum_{r \ge 1} \hat f(r) c_r(n), \phantom{mm}
g(n) = \sum_{s \ge 1} \hat g(s)c_s(n).
$$
Further suppose that there exists  $\alpha > 4$ such that
$$\big|\hat f(r)\big|, \big|\hat{g}(r) \big| \ll \frac{1}{r \log^\alpha r}$$
and $h$ is a positive integer.
Then for a positive integer $N$, we have
 $$\sum_{n \le N} f(n) g(n+h) =
N \sum_{r \ge 1} \hat f(r) \hat g(r) c_r(h) +
O \left( \frac{N}{(\log N)^{\alpha-4}} \right).
 $$
\end{thm}

Now for the Ramanujan coefficients coming from finite Ramanujan expansions
we first observe the following.

\begin{lem}\label{l2}
Suppose that the Ramanujan coefficients $\hat f(r)$ coming from the finite
Ramanujan expansion of $f$ satisfy
$$
\big|\hat f(r)\big| \ll \frac{1}{r \log^{\alpha}r}
\ \ \text{for some} \ \ \alpha>1.
$$
Then we have
$$
\big|f'(r)\big| \ll \frac{1}{\log^{\alpha-1}r}.
$$
Similarly if we assume
$$
\big|f'(r)\big| \ll \frac{1}{\log^{\beta}r}
\ \ \text{for some} \ \ \beta>1,
$$
then we have
$$
\big|\hat f(r)\big| \ll \frac{1}{r \log^{\beta-1}r}.
$$
\end{lem}

The proof uses the conversion formulas \eqref{hat f} and \eqref{f'} and partial
summation formula. For a detailed proof see Section \ref{lemmas}. Next we
prove the following theorem.

\begin{thm}\label{finite3}
Let $N$ be a positive integer and $f$ and $g$ be two arithmetical functions
for which we want to estimate the shifted convolution sums
$$
\sum_{n \le N} f(n) g(n+h)
$$
for a positive integer $h$. Further suppose that
$$
\big|f'(d)\big|, \big|g'(d)\big| \ll \frac{1}{\log^{\beta}d}
\ \ \text{for some} \ \ \beta>2,
$$
where $f'(d),g'(d)$ are as in \eqref{f} and \eqref{g}. Then we have
$$
\sum_{n \le N} f(n) g(n+h) = N \sum_{r=1}^{\infty}\hat{f}(r)\hat{g}(r)c_r(h)
+ O_{\beta,h}\Big(\frac{N}{\log^{\beta-2} N} \Big).
$$
\end{thm}

As an immediate corollary of \thmref{finite3}, we can now derive
(using \lemref{l2}) the following theorem, which is an analogoue of 
\thmref{BS}, in the setting of finite Ramanujan expansions. 

\begin{thm}\label{finite2}
Let $N$ be a positive integer and $f$ and $g$ be two arithmetical functions
for which we want to estimate the shifted convolution sums
$$
\sum_{n \le N} f(n) g(n+h)
$$
for a positive integer $h$. Further suppose that
$$
\big|\hat f(r)\big|,\big|\hat g(r)\big| \ll \frac{1}{r \log^{\alpha}r}
\ \ \text{for some} \ \ \alpha>3,
$$
where $\hat f(r),\hat g(r)$ are as in \eqref{hat f} and \eqref{hat g}. Then we have
$$
\sum_{n \le N} f(n) g(n+h) = N \sum_{r=1}^{\infty}\hat{f}(r)\hat{g}(r)c_r(h)
+ O_{\alpha,h}\Big(\frac{N}{\log^{\alpha-3} N} \Big).
$$
\end{thm}

%Note that here also we improve the error term by a factor of $\log N$.

\section{Proofs of the lemmas}\label{lemmas}

\begin{proof}[\bf Proof of \lemref{l1}]
The lemma follows from the known identity
$$
\frac{1}{d} \sum_{a=1}^d \zeta_d^{an}= \begin{cases}
1 & \text{if } \ d|n,\\
0 & \text{otherwise},
\end{cases}
$$
where $\zeta_d$ denotes a primitive $d$-th root of unity.
We then partition the sum in the left hand side in terms of
$\gcd$ and using the definition we write
$$
\frac{1}{d} \sum_{a=1}^d \zeta_d^{an}
= \frac{1}{d} \sum_{r|d} \sum_{a=1 \atop (a,d)=r}^d \zeta_d^{an}
= \frac{1}{d} \sum_{r|d} c_{d/r} (n)
= \frac{1}{d} \sum_{r|d} c_r (n).
$$
\end{proof}

\begin{proof}[\bf Proof of \lemref{l2}]
We use the formula \eqref{f'} and write
$$
\big|f'(d)\big|=\left| d\sum_{j\le {N\over d}}\mu(j)\hat{f}(jd) \right|
\ll \sum_{j\le {N\over d}}{1\over {j\log^{\alpha}(jd)}}.
$$
Without loss of generality we take $d>1$ and break the above sum in two different
cases: one for $d\le \sqrt{N}$ and the other for $d>\sqrt{N}$.
When $d\le \sqrt{N}$, we have $d\le {N\over d}$. Thus we have
$$
\sum_{j\le N/d}{1\over {j\log^{\alpha}(jd)}}
\ll \sum_{j\le d}{1\over {j\log^{\alpha} d}} 
+\sum_{d<j\le {N\over d}}{1\over {j\log^{\alpha} j}}
\ll_{\alpha} {1\over {\log^{\alpha-1} d}} + {1\over {\log^{\alpha-1}(N/d)}} 
\ll_{\alpha} {1\over {\log^{\alpha-1} d}}.
$$
Next, if $d>\sqrt{N}$ then ${N\over d}<d$. Thus
$$
\sum_{j\le N/d}{1\over {j\log^{\alpha}(jd)}}
\le \sum_{j\le d}{1\over {j\log^{\alpha} d}} 
\ll {1\over {\log^{\alpha-1} d}}.
$$
This completes the proof of the first part.
For the second part we use the formula \eqref{hat f} and write
$$
\big| \hat f(r) \big| = \left| \sum_{r|d, d\le N} \frac{f'(d) }{d} \right|
\ll {1\over r}\sum_{n\le N/r}{1\over {n\log^{\beta}(rn)}}.
$$
Again we take $r>1$ and break the above sum in two different
cases: one for $r \le \sqrt{N}$ and the other for $r>\sqrt{N}$.
If $ r\le \sqrt{N}$, then $r \le {N\over r}$ and hence
$$ 
\sum_{n\le N/r}{1\over {n\log^{\beta}(rn)}}
\ll \sum_{n\le r}{1\over {n\log^{\beta} r}} 
+ \sum_{r<n\le {N\over r}}{1\over {n\log^{\beta} n}}
\ll_{\beta} {1\over {\log^{\beta-1} r}} + {1\over {\log^{\beta-1}(N/r)}} 
\ll_{\beta} {1\over {\log^{\beta-1} r}}.
$$
Further if $r>\sqrt{N}$, then ${N\over r}<r$ and thus
$$ 
\sum_{n\le N/r}{1\over {n\log^{\beta}(rn)}}
\ll \sum_{n\le r}{1\over {n\log^{\beta} r}} 
\ll {1\over {\log^{\beta-1} r}}.
$$
This completes the proof.
\end{proof}

\section{Proofs of the theorems}

\begin{proof}[\bf Proof of \thmref{finite}]
We start with
$$
\sum_{n \le N} f(n) g(n+h)=\sum_{n\le N} \sum_{d|n}f'(d) \sum_{q|n+h}g'(q)
=\sum_{d\le N}f'(d)\sum_{q\le N+h}g'(q)
\sum_{{n\le N}\atop {{n\equiv 0\bmod d}\atop {n+h\equiv 0\bmod q}}}1.
$$
Now gathering by gcd and changing variables $d,q$ we get
$$
\sum_{n \le N} f(n) g(n+h)
=\sum_{{l|h}\atop {b:=-{h\over l}}} \sum_{d\le {N\over l}}f'(ld)
\sum_{q\le {{N+h}\over l}}g'(lq)
\sum_{{m\le {N\over {ld}}}\atop {m\equiv \overline{d}b \bmod q}}1.
$$
Here and from now on $\overline{d}$ denotes an inverse of $d$ modulo $q$.
Next we split the summations with conditions $dq\le N/l$ and $dq>N/l$
(the $\ast$ in $q-$sums abbreviates $(q,d)=1$ hereafter): 
\begin{align*}
\sum_{n \le N} f(n) g(n+h) = & \sum_{{l|h}\atop {b:=-{h\over l}}}
\sum_{d\le {N\over l}}f'(ld) \starsum_{q\le {N\over {ld}}}g'(lq)
\sum_{{m\le {N\over {ld}}}\atop {m\equiv \overline{d}b \bmod q}}1\\
& + O \Bigg( \sum_{{l|h}\atop {b:=-{h\over l}}}{1\over {l^{2\delta}}}
\sum_{d\le {N\over l}}{1\over d^{\delta}}
\sum_{{N\over {ld}}<q\le {{N+h}\over l}}{1\over q^{\delta}}
\sum_{{m\le {N\over {ld}}}\atop {md\equiv b \bmod q}}1 \Bigg),
\end{align*}
where the last sum is (thanks to condition $dq>N/l$)
$$
\ll_{\delta} \sum_{l|h}{1\over {l^{2\delta}}}\Big({l\over N}\Big)^{\delta}
\sum_{n\le {N\over l}}d(n)d(n+h/l)
\ll_{\delta,h} N^{1-\delta}(\log N)^2. 
$$
Here we have used the equivalent form of our hypothesis as per
\eqref{hypo-1} and \eqref{hypo-2}.
We now have the term 
$$
\sum_{{l|h}\atop {b:=-{h\over l}}}\sum_{d\le {N\over l}}f'(ld)
\starsum_{q\le {N\over {ld}}}g'(lq)
\sum_{{m\le {N\over {ld}}}\atop {m\equiv \overline{d}b \bmod q}}1
=\sum_{l|h}\sum_{d\le {N\over l}}f'(ld) \starsum_{q\le {N\over {ld}}}g'(lq)
\Big({N\over {ldq}}+O(1)\Big), 
$$
where the part with $O(1)$ contributes 
$$
\ll_{\delta} \sum_{l|h}{1\over {l^{2\delta}}}\sum_{d\le {N\over l}}{1\over {d^\delta}}\sum_{q\le {N\over {ld}}}{1\over {q^\delta}} 
\ll_{\delta} \sum_{l|h}{1\over {l^{2\delta}}}\sum_{d\le {N\over l}}{1\over {d^\delta}}\Big({N\over {ld}}\Big)^{1-\delta}
\ll_{\delta} N^{1-\delta}(\log N)
$$
if $\delta \neq 1$ and $\ll \log^2 N$ if $\delta=1$.
The main term is then coming from
$$
N\sum_{l|h}{1\over l}\sum_{d\le {N\over l}}{{f'(ld)}\over d}
\starsum_{q\le {N\over {ld}}}{{g'(lq)}\over q}
$$
which is written as
\begin{align*}
& N \sum_{l|h}{1\over l}\sum_{d\le {N\over l}}{{f'(ld)}\over d}
 \starsum_{q\le {{N+h}\over l}}{{g'(lq)}\over q}
 +O_{\delta}\Big( N \sum_{l|h}{1\over l^{1+2\delta}}
 \sum_{d\le {N\over l}}{1\over {d^{1+\delta}}}
 \sum_{{N\over {ld}}<q\le {{N+h}\over l}}{1\over {q^{1+\delta}}}\Big)\\
 &= N \sum_{l|h}l\sum_{d}{{f'(ld)}\over ld}
 \starsum_{q}{{g'(lq)}\over lq}
 +O_{\delta}\Big(N^{1-\delta}\sum_{l|h}{1\over l^{1+\delta}}
 \sum_{d\le {N\over l}}{1\over d}\Big)\\
 &= N \sum_{l|h}l\sum_{d}{{f'(ld)}\over ld}
\sum_{{q}\atop {(q,d)=1}}{{g'(lq)}\over lq}
 +O_{\delta}\Big(N^{1-\delta} \log N\Big).
\end{align*}
Next we use the following fundamental property of the M\"obius function
(See page 3 of \cite{CM}).
\begin{lem}\label{Mobius}
$$
\sum_{d|n} \mu(d) =
\begin{cases}
1 & \text{if } \ n=1,\\
0 & \text{otherwise}.
\end{cases}
$$
\end{lem}

Using \lemref{Mobius}, for $n=(q,d)$, we write
\begin{align*}
\sum_{l|h}l\sum_{d}{{f'(ld)}\over ld}
\sum_{{q}\atop {(q,d)=1}}{{g'(lq)}\over lq}
&= \sum_{l|h}l\sum_{t}\mu(t)\sum_{d'}{{f'(ltd')}\over {ltd'}}
\sum_{q'}{{g'(ltq')}\over {ltq'}}\\
&=\sum_{l|h}l \sum_{t}\mu(t)\hat{f}(lt)\hat{g}(lt)\\
&=\sum_{r=1}^{\infty}\hat{f}(r)\hat{g}(r)\sum_{{l|h}\atop {l|r}}l
\mu\left({r\over l}\right) \\
&=\sum_{r=1}^{\infty}\hat{f}(r)\hat{g}(r)c_r(h).
\end{align*}

Thus we get
$$
\sum_{n \le N} f(n) g(n+h) = N \sum_{r=1}^{\infty}\hat{f}(r)\hat{g}(r)c_r(h)
+ O_{\delta,h}\Big(N^{1-\delta} \log^2 N+1\Big).
$$

\end{proof}

\begin{rmk}
\rm
One can also follow steps of \cite{MS} and \cite{CMS}. The proof obtained
in this way will be a little shorter. However, this will only prove a weaker version
of this result. We briefly sketch it below.

Keeping the principle of \cite{MS} and \cite{CMS} in mind, we consider a
parameter $U$ tending to infinity which is to be chosen later. Then we write
$$
\sum_{n \le N} f(n) g(n+h) = A+B,
$$
where
$$ 
A:= \sum_{n \le N} \sum_{\substack{r,s\\ rs \le U}} \hat f(r) \hat g(s) c_r(n)c_s(n+h)
~\mbox{ and }~
B:=\sum_{n \le N} \sum_{\substack{r,s\\ rs > U}} \hat f(r) \hat g(s) c_r(n)c_s(n+h).
$$
As per our derivation in \cite{CMS}, we have
$$
A
=\begin{cases}
N \sum_{r \ge 1} \hat f(r) \hat g(r) c_r(h) +
O_h\left(\frac{N}{U^{1/2+\delta}}\right) + O(U^{1-\delta} \log^2 U)
& ~ \mbox{ if } ~ \delta < 1,\\
N \sum_{r \ge 1} \hat f(r) \hat g(r) c_r(h) +
O_h\left(\frac{N}{U^{3/2}}\right) + O(\log^3 U)
& ~ \mbox{ if } ~ \delta = 1,\\
N \sum_{r \ge 1} \hat f(r) \hat g(r) c_r(h) +
O_h\left(\frac{N}{U^{1/2+\delta}}\right) + O(1)
& ~ \mbox{ if } ~ \delta > 1.
\end{cases}
$$

Using \eqref{crn-exp1} and the hypotheses on $\hat f(r)$ and $\hat g(r)$ we write
$$
|B| \ll \sum_{\substack{r \le N ,s \le N+h\\ rs > U}}\frac{1}{(rs)^{1+\delta}}
\sum_{r' | r} r' \sum_{s' | s} s' 
\sum_{\substack{ n \le N \\ r'|n, s'|n+h}} 1.
$$
Next we put $r= r' n_r$ and $s=s' n_s$. Hence
\begin{align*}
|B| & \ll \sum_{r' \le N ,s' \le N+h} \frac{1}{(r's')^{\delta}}
\sum_{\substack{n_r ,n_s \\ n_r n_s > U/r's'}} \frac{1}{(n_r n_s)^{1+\delta}}
\sum_{\substack{ n \le N \\ r'|n,s'|n+h}} 1\\
& = \sum_{r' \le N ,s' \le N+h} \frac{1}{(r's')^{\delta}}
\sum_{t > U/r's'} \frac{d(t)}{t^{1+\delta}}
\sum_{\substack{ n \le N \\ r'|n, s'|n+h}} 1\\
& \ll  \sum_{r' \le N ,s' \le N+h} \frac{1}{(r's')^{\delta}}
\frac{\log (U/r's')}{(U/r's')^\delta}
\sum_{\substack{ n \le N \\ r'|n, s'|n+h}} 1\\
& \ll \frac{\log(UN^2)}{U^\delta} \sum_{n \le N} d(n) d(n+h)\\
& \ll_h \frac{N \log^2 N \log(UN^2)}{U^\delta}.
\end{align*}
To optimize the error terms, we choose
$$
U=\begin{cases}
N \log N & ~\mbox{ if }~\delta < 1,\\
N & ~\mbox{ if }~\delta = 1,\\
N^{1/\delta} (\log N)^{3/\delta} & ~\mbox{ if }~\delta > 1.
\end{cases}
$$
These choices yield
$$
 \sum_{n \le N} f(n) g(n+h) =
 \begin{cases}
N \sum_{r \ge 1} \hat f(r) \hat g(r) c_r(h) +
O(N^{1-\delta} (\log N)^{3-\delta}) & ~ \mbox{ if } ~ \delta < 1,\\
N \sum_{r \ge 1} \hat f(r) \hat g(r) c_r(h) +
O(\log^3 N) & ~ \mbox{ if } ~ \delta = 1,\\
N \sum_{r \ge 1} \hat f(r) \hat g(r) c_r(h)  +
O(1) & ~ \mbox{ if } ~ \delta > 1.
 \end{cases}
 $$
\end{rmk}

\begin{proof}[\bf Proof of \thmref{finite3}]
This proof starts off similarly. We just rewrite the hypotheses on $f'$ and $g'$ as
$$
\big| f'(d) \big|, \big| g'(d) \big| \ll {1\over {1+\log^\beta d}}
$$
and obtain
\begin{align*}
\sum_{n \le N} f(n) g(n+h) = & \sum_{{l|h}\atop {b:=-{h\over l}}}
\sum_{d\le {N\over l}}f'(ld) \starsum_{q\le {N\over {ld}}}g'(lq)
\sum_{{m\le {N\over {ld}}}\atop {m\equiv \overline{d}b \bmod q}}1\\
& + O \Bigg( \sum_{{l|h}\atop {b:=-{h\over l}}}
\sum_{d\le {N\over l}}{1\over {1+\log^\beta(ld)}}
\starsum_{{N\over {ld}}<q\le {{N+h}\over l}}{1\over {1+\log^\beta(lq)}}
\sum_{{m\le {N\over {ld}}}\atop {m\equiv \overline{d}b \bmod q}}1 \Bigg)\\
= & \sum_{{l|h}\atop {b:=-{h\over l}}}
\sum_{d\le {N\over l}}f'(ld) \starsum_{q\le {N\over {ld}}}g'(lq)
\sum_{{m\le {N\over {ld}}}\atop {m\equiv \overline{d}b \bmod q}}1 + O(R_1),
\end{align*}
where
$$
R_1:=
\sum_{{l|h}\atop {b:=-{h\over l}}}
\sum_{d\le {N\over l}}{1\over {1+\log^\beta(ld)}}
\starsum_{{N\over {ld}}<q\le {{N+h}\over l}}{1\over {1+\log^\beta(lq)}}
\sum_{{m\le {N\over {ld}}}\atop {m\equiv \overline{d}b \bmod q}}1.
$$
To estimate $R_1$ we separate according to the cases $d \le \sqrt{N}$
and $d > \sqrt{N}$. If $d \le \sqrt{N}$ then
$$
\log^{\beta}(lq)\gg_{\beta} \log^{\beta} N
$$
for all $q>N/ld$, while $d > \sqrt{N}$ implies
$$
\log^{\beta}(ld)\gg_{\beta} \log^{\beta} N.
$$
This yields
\begin{align*}
R_1 & \ll_{\beta} \sum_{{l|h}\atop {b:=-{h\over l}}} \frac{1}{\log^\beta N}
\sum_{d\le {N\over l}}\sum_{{N\over {ld}}<q\le {{N+h}\over l}}
\sum_{{m\le {N\over {ld}}}\atop {md\equiv b \bmod q}}1\\
& \ll_{\beta} {1\over {\log^\beta N}}\sum_{l|h}\sum_{n\le {N\over l}} d(n) d(n+h/l)\\
& \ll_{\beta,h} {1\over {\log^{\beta} N}} \sum_{l|h} {N\over l} \log^2 {N\over l}\\
& \ll_{\beta,h} {N\over {\log^{\beta-2} N}}.
\end{align*}
Here we used the asymptotic estimate
$$
\sum_{n \le N/l} d(n) d(n+h/l) \sim
\frac{6}{\pi^2} \sigma_{-1}(h/l) \frac{N}{l} \log^2(N/l),
$$
due to Ingham \cite{AEI}. So now we are left to estimate
$$
\sum_{{l|h}\atop {b:=-{h\over l}}}
\sum_{d\le {N\over l}}f'(ld) \starsum_{q\le {N\over {ld}}}g'(lq)
\sum_{{m\le {N\over {ld}}}\atop {m\equiv \overline{d}b \bmod q}}1
= N\sum_{l|h}{1\over l}\sum_{d\le {N\over l}}{{f'(ld)}\over d}
\starsum_{{q\le {N\over {ld}}}}{{g'(lq)}\over q} 
 +O(R_2),
$$
where
$$
R_2:=\sum_{l|h}\sum_{d\le {N\over l}}|f'(ld)|\sum_{q\le {N\over {ld}}}|g'(lq)|.
$$
Here we used the fact that
$$
\sum_{{m\le {N\over {ld}}}\atop {m\equiv \overline{d}b\bmod q}}1 = {N\over {ldq}}+O(1).
$$
Using the hypothesis we get that
$$
R_2 \ll \sum_{l|h}\sum_{d\le {N\over l}}{1\over {1+\log^\beta(ld)}}
\sum_{q\le {N\over {ld}}}{1\over {1+\log^\beta(lq)}}.
$$
To treat the $q$-sum on the right hand side we split as follows:
$$
\sum_{q\le {N\over {ld}}}{1\over {1+\log^\beta(lq)}}
\ll_{\beta} \sum_{q\le {1\over l}\sqrt{N\over d}}1
+{1\over {1+\log^{\beta}(N/d)}}\sum_{{1\over l}\sqrt{N\over d}<q\le {N\over {ld}}}1,
$$
where we used that
$$
q>{1\over l}\sqrt{N\over d} \Rightarrow
{1\over {1+\log^\beta(lq)}}\ll_{\beta} {1\over {1+\log^{\beta}(N/d)}}.
$$
Hence
\begin{align*}
R_2 & \ll_{\beta} \sum_{l|h}\sum_{d\le {N\over l}}{1\over {1+\log^\beta(ld)}}
\ {N\over {ld(1+\log^{\beta}(N/d))}}\\
& \ll_{\beta} {N\over {\log^\beta N}}\sum_{l|h}{1\over l}
\left(1+\sum_{1<d\le \sqrt{N}}{1\over {d\log^\beta(ld)}}
+\sum_{\sqrt{N}<d\le {N\over l}}{1\over {d(1+\log^{\beta}(N/d))}}\right)\\
& \ll_\beta {N\over {\log^\beta N}}\sum_{l|h}{1\over l} \sum_{d\le {N\over l}}{1\over d}\\
& \ll_{\beta,h} \frac{N}{\log^{\beta-1} N}.
\end{align*}
Thus so far we have obtained that
$$
\sum_{n \le N} f(n) g(n+h) = N\sum_{l|h}{1\over l}\sum_{d\le {N\over l}}{{f'(ld)}\over d}
\starsum_{{q\le {N\over {ld}}}}{{g'(lq)}\over q} +
O_{\beta,h}\left( \frac{N}{\log^{\beta-2} N} \right).
$$
Now we essentially repeat what we did in the proof of \thmref{finite} and write
$$
N\sum_{l|h}{1\over l}\sum_{d\le {N\over l}}{{f'(ld)}\over d}
\starsum_{{q\le {N\over {ld}}}}{{g'(lq)}\over q}
= M + O(R_3),
$$
where
$$
M:=N\sum_{l|h}{1\over l}\sum_{d\le {N\over l}}{{f'(ld)}\over d}
\starsum_{{q\le {{N+h}\over l}}}{{g'(lq)}\over q}
=N\sum_{r = 1}^\infty\hat{f}(r)\hat{g}(r)c_r(h)
$$
and
$$
R_3:=N\sum_{l|h}{1\over l}\sum_{d\le {N\over l}}{{|f'(ld)|}\over d}
\sum_{{N\over {ld}}<q\le {{N+h}\over l}}{{|g'(lq)|}\over q}.
$$
Using the hypothesis we get that
\begin{align*}
R_3
& \ll N\sum_{l|h}{1\over l}\sum_{d\le {N\over l}}{1\over d(1+\log^\beta(ld))}
\sum_{{N\over {ld}}<q\le {{N+h}\over l}}{1\over q(1+\log^\beta(lq))}\\
& \ll N\sum_{l|h}{1\over l}\sum_{d\le {N\over l}}{1\over d(1+\log^\beta(ld)) (1+\log^\beta(N/d))}
\sum_{{N\over {ld}}<q\le {{N+h}\over l}}{1\over q}\\
& \ll \frac{N}{\log^\beta N} \sum_{l|h}{1\over l}\sum_{d\le {N\over l}}{1\over d}
\sum_{{N\over {ld}}<q\le {{N+h}\over l}}{1\over q}\\
& \ll_h \frac{N}{\log^{\beta-2} N}.
\end{align*}
This completes the proof.
\end{proof}

\section{Concluding remarks}

The method outlined in this paper will undoubtedly have
further applications as the theory moves forward.  It offers
us yet another way to approach these general convolution sums.
The technical issues regarding absolute convergence of infinite
series that complicated our earlier work have now been simplified
through the use of finite Ramanujan expansions. 

As demonstrated in the paper, {\it all} arithmetical functions 
now afford a {\it finite} Ramanujan expansion. Convolution 
sums lie at the heart of analytic number theory. Earlier, 
there have been attempts to study such sums. Our paper 
offers yet another route to this study. We expect to
investigate in future work further refinements of the theory.\\

{\bf Acknowledgements:} We would like to thank the referee for
useful remarks which improved the presentation of this article.
We would also like to thank the organisers of conferences
`Leuca 2016' and `CNTA 2016' for their kind hospitality,
where the final part of this work was done.

\end{document}